\documentstyle[12pt,twoside]{amsart}

\newtheorem{thm}{Theorem}[section]

\newtheorem{lem}[thm]{Lemma}

\newtheorem{prop}[thm]{Proposition}

\theoremstyle{definition}
\newtheorem{defn}[thm]{Definition}
\newtheorem{condition}[thm]{Condition}
\newtheorem{say}[thm]{}
\newtheorem{exmp}[thm]{Example}

\newtheorem{rem}[thm]{Remark}          
\newtheorem{ack}{Acknowledgments}        
\newtheorem{notation}[thm]{Notation}

\theoremstyle{remark}

\renewcommand{\c}[0]{{\mathbb C}}  

\renewcommand{\o}[0]{{\cal O}} 
\newcommand{\z}[0]{{\mathbb Z}}

\renewcommand{\r}[0]{{\mathbb R}} 

\renewcommand{\a}[0]{{\mathbb A}}

\newcommand{\p}[0]{{\mathbb P}}

\newcommand{\q}[0]{{\mathbb Q}}

\newcommand{\qtq}[1]{\quad\mbox{#1}\quad}

\newcommand{\pic}[0]{\operatorname{Pic}}

\newcommand{\mult}[0]{\operatorname{mult}}

\newcommand{\supp}[0]{\operatorname{Supp}}    
\newcommand{\red}[0]{\operatorname{red}}

\newcommand{\sing}[0]{\operatorname{Sing}}

\newsymbol\subsetneq 2328
\newsymbol\onto 1310
\def\into{\DOTSB\lhook\joinrel\rightarrow}
\newsymbol\twoheadrightarrow 1310

\begin{document}
\bibliographystyle{amsplain}

\title{Real Algebraic Threefolds IV.\\ Del Pezzo Fibrations}
\author{J\'anos Koll\'ar}

\maketitle
\tableofcontents

\section{Introduction}

This paper continues the study of the topology of real algebraic
threefolds begun in \cite{rat1, rat2, rat3}, but the current work is  
independent of the previous ones in its methodology.

The present aim is to understand the topology of  the set of real
points of threefolds which admit a morphism to a curve whose general
fiber is a rational surface.   This class of threefolds also appears as
one of the 4 possible outcomes of the minimal model program (cf.\
\cite[Sec.\ 3.7]{KoMo98}). 
 Our main theorem gives a nearly complete description of the possible
topological types of the set of real points of such a threefold in the
orientable case.

\begin{thm}\label{main.thm} Let $X$ be a smooth projective real
algebraic threefold such that the set of real points $X(\r)$ is
orientable.
  Assume that  there is a  morphism  $f:X\to C$  onto a real algebraic
curve $C$ whose 
 general fibers   are rational surfaces. Let $M\subset X(\r)$ be any
connected component. Then 
$$ M\sim N\ \#\ a\r\p^3 \ \#\ b(S^1\times S^2)\qtq{for some $a,b\geq
0$,}
$$ where   one of the following holds.
\begin{enumerate}
\item $N$ is a connected sum of lens spaces (cf.\ (\ref{lens.defn})).
\item  $N$ is Seifert fibered over a topological surface (cf.\
(\ref{seif.defn})).
\item   $N$ is either
 an $S^1\times S^1$-bundle over $S^1$ or  is doubly covered by such a
bundle. 
\end{enumerate}
\end{thm}

\begin{say}[First reduction step]\label{red.step.say}{\ }

As with many results for 3--folds, the proof starts with a  suitable
minimal model program. The  minimal model program over $\r$  was studied
in greater generality  in
\cite{rat2} so I just outline the main conclusions. The end result is
that (\ref{main.thm}) follows from  (\ref{main.dp.thm}) and for the rest
of the paper we do  not have  to know anything about the 3--dimensional
minimal model program.

Let us run the relative minimal model program for $X\to C$ over $\r$
(cf.\ \cite[Sec.\ 3.7]{KoMo98} or \cite[Sec.\ 3]{rat2}). The change in
the topology of real points in the course of the program is described in
\cite[1.2]{rat2}.   There are two possibilities for the final step of
the program.

First, we may get   a conic bundle over a surface. These   are
described in \cite{rat3} and we obtain cases 1 or 2 of
(\ref{main.thm}).

Otherwise the program ends with a morphism $f^m:X^m\to C$ which has the
following properties:
\begin{enumerate}
\item $X^m$ has only isolated singularities which are $\q$-factorial
over $\r$ (that is, if $D$ is a Weil divisor defined over $\r$ then
$mD$ is Cartier for some $m>0$). 
\item $-K_{X^m}$ is $f^m$-ample.
\item $K_{X^m}$  is Cartier at all real points.
\item Every fiber of $f$ is irreducible (over $\r$).
\item The topological normalization $\overline{X^m(\r)}$  of $X^m(\r)$ 
is a 3--manifold.
\end{enumerate} In view of \cite[1.2]{rat2} and \cite[1.1]{rat3},
(\ref{main.thm}) is reduced to proving the following:
\end{say}

\begin{thm}\label{main.dp.thm} Let $X$ be a  projective real algebraic
threefold  and $f:X\to C$  a morphism onto a real algebraic curve $C$.
Assume that  $X$ satisfies the conditions (\ref{red.step.say}.1--5) and
that $\overline{X(\r)}$ is orientable. Let $M\subset \overline{X(\r)}$
be any connected component. Then 
$$ M\sim N\ \#\ a\r\p^3 \ \#\ b(S^1\times S^2)\qtq{for some $a,b\geq
0$,}
$$ where   one of the following holds.
\begin{enumerate}
\item $N$ is a connected sum of lens spaces.
\item   $N$ is either
 an $S^1\times S^1$-bundle over $S^1$ or  is doubly covered by such a
bundle. 
\end{enumerate}
\end{thm}

\begin{rem} The conclusion of the theorem can probably be
considererably improved. 

In case (\ref{main.dp.thm}.1) I do not have any examples where $N$ is
not a single lens space. A more detailed analysis of the present
methods may prove that this is always the case.

$(S^1\times S^1)$-bundles over $S^1$ are classified by the monodromy map
on $H_1(S^1\times S^1,\z)\cong \z^2$.  If this map is hyperbolic then 
$N$ has a geometry modelled on {\it Sol} (cf.\ \cite[p.470]{Scott83}).
I believe that  $N$ can never be of this type but the methods of my
proof do not say anything about the monodromy.
\end{rem}

\begin{say}[Second reduction step]\label{2nd.red.say}{\ }

In this step we reduce the proof of (\ref{main.dp.thm}) to the study of
the singular fibers of $f$. 

Let $A\sim S^1$ be a connected component of $C(\r)$ and  
$p_1,\dots,p_s\in A$   the points (in cyclic order) over which $f$ is
not smooth. For each $i$ pick a point $q_i\in (p_i,p_{i+1})$. Then
$X_{q_i}:=f^{-1}(q_i)$ is a smooth Del Pezzo surface  and $X_{q_i}(\r)$
is orientable. Thus 
$X_{q_i}(\r)$ is either $S^1\times S^1$ or a disjoint union of copies
of $S^2$
by a result of  \cite{Comessatti14}. (See also
\cite[V.3.4,VI.4.6 and VI.6.3]{Silhol89}.)
 Gluing 3-manifolds along such surfaces is  a relatively
simple operation, thus one can expect to  get a good description of
$X(\r)$ by describing the pieces
$Z_i:=(f^{-1}[q_{i-1},q_i])(\r)$ for every $i$; see
(\ref{main.thm.pf.2}). 

$f:Z_i\to [q_{i-1},q_i]$ is a function whose only   critical value is
$p_i$. Thus $Z_i$ can be viewed as a regular neighborhood of the
critical level set $(f^{-1}(p_i))(\r)$. In fact, once we have a
topological  description of    $(f^{-1}(p_i))(\r)$, it is easy to
figure out what
$Z_i$ is, at least up to finite ambiguity. This is done in 
(\ref{main.thm.pf.1}). 

The complex projective surface $S_i:= f^{-1}(p_i)$ is a ``singular Del
Pezzo" surface which appears as a degeneration of smooth Del Pezzo
surfaces. Quite a lot is known about such surfaces. A structural theory
of   normal  Del Pezzo surfaces is developed in \cite{Keel-Mc}.
Degenerations of
$\p^2$ and some other  Del Pezzo surfaces are studied in
\cite{Manetti1, Manetti2}. Unfortunately, the conclusion of these works
is that a complete classification of such singular Del Pezzo surfaces
is not   feasible because of the combinatorial complexity of the
problem. On the other hand, the methods developed by these and other
authors can be used to answer  many questions about singular Del Pezzo
surfaces.

The main part of this work is devoted to understanding the topology of
the set of real points of certain singular Del Pezzo surfaces. It
should be emphasized that, as opposed to almost all previous studies, we
have to consider  nonnormal  surfaces as well. In fact, as far as the
topology of the real points is concerned, normal surfaces present no
difficulties. 

Nonnormal   Del Pezzo surfaces whose canonical divisor is Cartier have
been enumerated in \cite{Reid94}.  Many of the irreducible ones appear
as singular fibers in our case. For all such examples one can  perform
a small   perturbation  of $f:Z_i\to [q_{i-1},q_i]$ such  that we
obtain a Morse function whose fibers stay {\it real algebraic} Del
Pezzo surfaces. These are again very easy to understand topologically. 

I have been unable to find any other nonnormal irreducible examples,
thus it is possible that  the work of  sections 2--4 is entirely
superfluous. On the other hand, the methods of these sections lead to
a   description of the singular fibers in many cases and it may be
possible to develop them further to obtain a complete list of the
 irreducible  fibers occurring in (\ref{main.dp.thm}).

Some interesting examples are given in section 7.
\end{say}

\begin{say}[The fundamental group of a real Del Pezzo surface]{\ }

Let $S:=f^{-1}(p_i)$ be one of the fibers and
$Z:=(f^{-1}[q_{i-1},q_i])(\r)$ the corresponding piece of the above
decomposition of $X(\r)$. 

Here I would like to outline a simple argument giving some information
about the  fundamental group of $Z$. Since $S(\r)$ is a retract of
$Z$, we see that $\pi_1(Z)\cong \pi_1(S(\r))$.

Let $B\subset S$ be the singular locus and   $U_j$  the connected
components of $S(\r)\setminus B(\r)$.  
$\pi_1(B(\r))$ is a free group and each $U_j$ gives one relation. Hence
we see that the number of relations of $\pi_1(S(\r))$ is bounded by the
number of  connected components of $S(\r)\setminus B(\r)$.

Since $S$ is a limit of smooth Del Pezzo surfaces, $|-K_S|$ has a
member; call it $D$. (This is a bit nebulous since $S$ may not even be
normal; see (\ref{degen.of.dp.say}) for its precise meaning.) Let
$h:S'\to S$ be the minimal resolution of $S$. It is easy to see that
$|-K_{S'}|$ has a unique member $D'$ such that $h_*(D')=B+D$.  Let
$f:S'\to S''$ be a minimal model and set $D'':=f_*D'$. 

$f$ is a composite of blow ups and one checks that
$$ |\pi_0(S''(\r)\setminus D''(\r))|
\geq  |\pi_0(S'(\r)\setminus D'(\r))|.
$$ Using \cite[3.39]{KoMo98} one easily checks that every fiber of $h$
is either contained in $\supp D'$ or is disjoint from it. This gives
that
$$ |\pi_0(S'(\r)\setminus D'(\r))|
\geq  |\pi_0(S(\r)\setminus (B+D)(\r))|.
$$ Thus we are led to counting  the  number of  connected components of
$S''(\r)\setminus D''(\r)$. A typical case is when $S''=\p^2$ and $D''$
is a cubic curve. We see that there are at most 4 connected components,
and 4 is achieved when $D''$ consists of 3  lines.   If we take into
account that we want to use only 
$S(\r)\setminus B(\r)$ and not $S(\r)\setminus (B+D)(\r)$, we  obtain
that  $S(\r)\setminus B(\r)$  has at most 2 connected components in
this case.

The other possible minimal models $S''$ can be similarly treated. (Ad
hoc arguments are needed to exclude the case when $S''$ is a minimal
ruled surface with negative section $E$ and $D''=2E+(\mbox{many
fibers})$.) At the end we obtain that
\begin{enumerate}
\item[($*$)] $\pi_1(Z)$ is the free product of groups with 1 relation.
\end{enumerate} On the one hand ($*$) is quite strong since most
3--manifold groups do not have this property. 
On the other hand, ($*$) is
not strong enough to exclude all hyperbolic 3--manifolds. 

The proof of (\ref{main.dp.thm}) given in sections 2--6 is essentially
an elaboration of this approach. The steps of going from $S$ to $S'$ and
$S''$ are studied in more detail. At the end we obtain a rather complete
geometric description of $S(\r)$.
\end{say}

\begin{rem}[PL three manifolds]\label{1.basic.top.facts} In this paper
I usually work with {\it piecewise linear  manifolds} (\cite{RoSa82} is
a  good introduction). Every real algebraic variety carries a natural
PL structure (cf.\ \cite[Sec.9.2]{BCR87}).

 In dimension 3 every compact topological 3--manifold carries a unique
PL--manifold structure (cf.\ \cite[Sec.\ 36]{Moise77}) and  a
PL--structure behaves very much like a differentiable  structure. For
instance, let $M^3$ be a PL 3--manifold,
$N$ a compact PL--manifold of  dimension 1 or 2 and $g:N\into M$ a
PL--embedding. Then a suitable open  neighborhood of $g(N)$ is
PL--homeomorphic to a real vector bundle over $N$  (cf.\ \cite[Secs.\
24 and 26]{Moise77}).  (The technical definition of  such neighborhoods
is given by the notion of  {\it regular neighborhood}, see
\cite[Chap.3]{RoSa82}).  If $f:M\to N$ is a PL--map and $X\subset N$ a
compact subcomplex then there is a regular neigborhood $X\subset
U\subset N$ such that
$f^{-1}(U)$  is a regular neigborhood $f^{-1}(X)\subset M$ 
 (cf.\ \cite[2.14]{RoSa82}). 
\end{rem}

\begin{defn}[Lens spaces]\label{lens.defn}{\ } For relatively prime
$0<q<p$ consider the action
$(x,y)\mapsto (e^{2\pi i/p}x, e^{2\pi iq/p}y)$ on the unit sphere
$S^3\sim (|x^2|+|y^2|=1)\subset \c^2$. The quotient is a 3--manifold
called the {\it lens space}  $L_{p,q}$. 

Another way to obtain lens spaces is to glue two solid tori together.
The result is   a lens space,  $S^3$ or $S^1\times S^2$.  Sometimes one
writes $L_{1,0}=S^3$ and $L_{0,1}=S^1\times S^2$, though these are
usually not considered lens spaces. (See, for instance,
\cite[p.20]{Hempel76}.)
\end{defn}

\begin{defn}[Seifert fiber spaces]\label{seif.defn}{\ } A 3--manifold
$M$ is called {\it Seifert fibered} if there is a morphism  
$f:M\to F$ to a topological surface such that every $P\in F$ has a
neighborhood  $P\in U\subset F$ such that $f:f^{-1}(U)\to U$ is fiber
preserving homeomorphic to  one of the normal forms $f_{c,d}$ defined
below.

Let $S^1\subset \c$ be the unit circle with coordinate $u$ and
$D^2\subset \c$ the closed unit disc with coordinate $z$. For a pair of
integers $c,d$ satisfying  $0\leq c<d$ and
$(c,d)=1$, define
$$ f_{c,d}:S^1\times D^2\to D^2\qtq{by}
 f_{c,d}(u,z)=u^cz^d.
$$
$f_{c,d}$ restricts to a fiber bundle
$S^1\times (D^2\setminus \{0\})\to D^2\setminus \{0\}$. The fiber of
$f_{c,d}$ over the origin is still $S^1$, but 
$f_{c,d}^{-1}(0)$ has multiplicity $d$. 

(This is the classical definition of Seifert fibered spaces, which is
slightly more restrictive than the one in \cite{Scott83}.)
\end{defn}

\begin{ack}  I   thank  M.\ Kapovich  and S.\ Kahrlamov for answering
my numerous questions about  3-manifold topology and real algebraic
geometry.   Partial financial support was provided by  the NSF under
grant number  DMS-9622394. 
\end{ack}

\section{Real Surfaces with Du Val Singularities}

The minimal model theory of real surfaces has been studied in detail in
the papers \cite{Comessatti14, Silhol89, ras}. It is not difficult to
generalize these results to the case when we allow Du Val singularities
(\ref{dv.defn}). 
 These results were explained in
\cite{rat3}.

\begin{defn}\label{dv.defn} Let $(0\in S)$ be a normal surface
singularity over $\r$ with minimal resolution $g:S'\to S$. 
  $(0\in S)$ is called a {\it    Du Val} singularity  (or rational
double point) iff $K_{S'}\cong g^*K_S$. Equivalently, 
$(0\in S)$ is      Du Val iff every
$g$-exceptional curve is a smooth rational curve with selfintersection
$-2$.  (See \cite{Reid85} or \cite[4.2]{KoMo98} for the relevant
background on Du Val singularities over $\c$.)

It is not hard to see that every  real Du Val singularity is real
analytically equivalent to one of the following normal forms (cf.
\cite[I.17.1]{AGV85} or \cite[Sec.\ 4.2]{KoMo98}).
\begin{description}
\item[$A_n^+$] $(x^2+y^2-z^{n+1}=0)$ for $n\geq 2$,
\item[$A_n^-$] $(x^2-y^2-z^{n+1}=0)$  for $n\geq 0$,
\item[$A_n^{++}$] $(x^2+y^2+z^{n+1}=0)$ for $n$ odd,
\item[$D_n^+$] $(x^2+y^2z+z^{n-1}=0)$ for $n\geq 4$,
\item[$D_n^-$] $(x^2+y^2z-z^{n-1}=0)$ for $n\geq 4$,
\item[$E_6^+$] $(x^2+y^3+z^4=0)$,
\item[$E_6^-$] $(x^2+y^3-z^4=0)$,
\item[$E_7$] $(x^2+y^3+yz^3=0)$,
\item[$E_8$] $(x^2+y^3+z^5=0)$.
\end{description}
\end{defn}

\begin{defn}\label{-1.curve.def}
 Let $S$ be a surface with Du Val singularities. A curve $C\subset S$
is called a {\it $(-1)$-curve} if its birational transform on the
minimal resolution is a smooth rational curve with selfintersection
$-1$. Thus $(C\cdot K_S)=-1$.
\end{defn}

\begin{defn}\label{(1,m)-bu.defn}
 Let $P\in S(\r)$ be a smooth real point and $x,y$ local coordinates at
$P$.  The surface
$S'\subset S\times \p^1_{(u:v)}$ given by equation $ux-vy^m=0$ is called
a  {\it $(1,m)$-blow up}  of $P$  on $S$.  For $m=1$ this is the
ordinary blow up. A $(1,m)$-blow up has a unique singular point of type
$A_{m-1}^-$ at $(0,0,0,1)$. 

It should be noted that for $m\geq 2$ the 
$(1,m)$-blow up does depend on the choice of the local coordinates. 

It is frequently better to think of a   $(1,m)$-blow up  as follows.
First blow up $0\in S$ to get $S_1\to S$.  Then blow up a point on the
exceptional divisor of $S_1\to S$ to obtain
$S_2\to S_1$. Then  blow up a point on the exceptional divisor of
$S_2\to S_1$  which is not on (the birational transform of) any
previous exceptional divisor to obtain
$S_3\to S_2$. After $m$-times  we have $m$ exceptional curves in the
following configuration:
$$
\stackrel{-1}{\circ} - \stackrel{-2}{\circ} - \dots  -
\stackrel{-2}{\circ}.
$$ We can now contract all the $(-2)$-curves    to get a 
$(1,m)$-blow up. This shows that the exceptional curve of a
$(1,m)$-blow up is a $(-1)$-curve.

If $P,\bar P\in S(\c)$ are smooth and conjugate complex points, then we
can choose conjugate coordinate systems to do a $(1,m)$-blow up at both
points. The result is again a real algebraic surface with a  conjugate
pair of $A_{m-1}$-points (for nonreal points the signs in the equations
do not matter). 
\end{defn}

\begin{defn}\label{F^r-def}
 Let $F_1,F_2$ be real algebraic surfaces with Du Val singularities and
assume that $g:F_1\to F_2$ is a composite of
$(1,m)$-blow ups. A $(1,m)$-blow up of a conjugate point pair is an
isomorphism in the neighborhood of the real points, so if a 
$(1,m)$-blow up of a conjugate point pair is followed by a
$(1,m)$-blow up of a real point then their order can be reversed.
Repeating if necessary, $g$  can be factored  uniquely as
$$ g:F_1\stackrel{g^c}{\to} F^r \stackrel{g^r}{\to} F_2
$$ where $g^c$ is a composite of
$(1,m)$-blow ups of  conjugate point pairs and 
$g^r$ is a composite of
$(1,m)$-blow ups of  real points.
\end{defn}

\begin{defn} Let $F$ be a  normal projective surface such that $K_F$ is
$\q$-Cartier. $F$ is called a (singular) {\it Del Pezzo} surface if
$-K_F$ is ample.
$F$ is called a {\it weak Del Pezzo} surface if $-K_F$ is nef and big.
A morphism   $g:F\to C$ is called a {\it conic bundle} if every fiber
is isomorphic to a plane conic (which can be smooth, a pair of
intersecting lines or a double line). Every conic bundle can be
embedded into a $\p^2$-budle over $C$ such that the fibers become
conics.
\end{defn}

  Combining the results of \cite[Sec.\ 9]{rat3} with (\ref{F^r-def}) we
obtain the following.

\begin{thm}\label{surf.dv.mmp}
 Let $F$ be a   projective  surface  over
$\r$  with Du Val singularities. Then there are surfaces and morphisms
$$ g:F\stackrel{g^c}{\to} F^r \stackrel{g^r}{\to} F^*
$$ with the following properties. 
\begin{enumerate}
\item $F^r$ and $F^*$ are  projective  surfaces  over
$\r$  with Du Val singularities.
\item  $g^c$ is a composite of
$(1,m)$-blow ups of  conjugate point pairs. In particular
$F(\r)\cong F^r(\r)$.
\item  $g^r$ is a composite of
$(1,m)$-blow ups of  real points.
\item $F^*$ falls in  one of the following 3 cases:
\begin{enumerate}
\item[(C)] (Conic bundle) $\rho(F^*)=2$ and
$F^*$ is a conic bundle over a smooth curve $A$.

\item[(D)] (Del Pezzo surface) $\rho(F^*)=1$ and $-K_{F^*}$ is ample.
\item[(N)] (Nef canonical class)  $K_{F^*}$ is nef.\qed
\end{enumerate}
\end{enumerate}
\end{thm}

Next we collect some auxiliary results that are needed elsewhere.

\begin{lem}\label{dp.pic.divis}
 Let $K$ be a field   and $F$ a Del Pezzo surface over $K$ with Du Val
singularities. Assume that
$\rho(F)=1$ and that there is an effective Cartier divisor
$0\neq B\subset F$ such that $-(K+B)$ is ample. 

Then $F$ is isomorphic to $\p^2$ or to a quadric hypersurface in $\p^3$.
\end{lem}

Proof. Let $H$ be a generator of $\pic(F)/(\mbox{torsion})$ and write 
$B\equiv bH$, $K_F\equiv -aH$.   This implies that $(K_F^2)=a^2(H^2)$.
By assumption $a>b>0$, so $a\geq 2$. Since
$(K_F^2)\leq 9$, thare three possibilities:
\begin{enumerate}
\item $a=3$ and $(H^2)=1$. Then $F_{\bar K}\cong \p^2$ and $F$ contains
a line defined over $K$, thus $F\cong \p^2$.
\item $a=2$ and $(H^2)=2$. Then $F_{\bar K}$ is a quadric and $\o(1)$
is defined over $K$.  Thus $F$ itself is a quadric.
\item  $a=2$ and $(H^2)=1$. Then 
$2p_a(B)-2=B\cdot (B+K_F)=H\cdot (-H)=-1$, which is impossible.\qed
\end{enumerate}

\begin{lem}\label{cb.with.sect}  Let $F$ be a real surface with Du Val
singularities and $\rho(F)=2$.  Let $F\to A$ be a conic bundle with a
section $H$. Then $H$ intersects  every singular fiber at a singular
point of $F$.
\end{lem}

Proof. A section can never intersect a multiple fiber at a smooth point
of $F$. If $f^{-1}(a)$ is a pair of conjugate lines then their
intersection point $P$ is the only real point, hence $H$ passes 
through $P$. If $P$ is a smooth point of $F$ then the intersection
number of
$H$ and $f^{-1}(a)$ is at least 2, which is impossible.\qed

\section{The Basic Set-up}

\begin{notation}\label{bas.setup.not} Let $X$ be a real algebraic
threefold and $f:X\to C$ a proper morphism to a smooth real algebraic
curve.  As a   generalization of the conditions of (\ref{main.dp.thm})
we assume the following:
\begin{enumerate}
\item $X$ has isolated $\q$-factorial  singularities,
\item $-K_X$ is $\q$-Cartier and $f$-ample, and
\item every fiber of $f$ is irreducible (over $\r$).
\end{enumerate}

Let $0\in C(\r)$ be a point and $X_0:=f^{-1}(0)$ the fiber  over $0$.
\end{notation}

\begin{lem}\label{3-fibers.lem}
 Notation as above. There are 3 possibilities for  $X_0$.
\begin{enumerate}
\item $X_0$ is reduced and geometrically  irreducible.
\item $X_0=mZ_0$ for some $m\geq 2$ where $Z_0$ is reduced and
geometrically  irreducible.
\item $X_0=m(Z_0+\bar Z_0)$ for some $m\geq 1$ where $Z_0$ and $\bar
Z_0$ are conjugate, reduced and irreducible.
\end{enumerate}
\end{lem}

Proof.  Write $(X_0)_{\c}=\sum a_iZ_i$ as the sum of its irreducible
and reduced components. For any $Z_i$ let $Z'_i:=Z_i$ if $Z_i$ is
defined over $\r$ and 
$Z'_i:=Z_i+\bar Z_i$ otherwise.
$Z'_i$ is defined over $\r$ hence  $(X_0)_{\c}=mZ'_i$ for some
$m\geq 1$ since $X_0$ is irreducible over $\r$.\qed

\begin{say}\label{degen.of.dp.say}   Let $X$ be a real algebraic
threefold and $f:X\to C$ a proper morphism to a smooth real algebraic
curve,
$0\in C(\r)$  a point and
$Z:=f^{-1}(0)$ the fiber  over $0$.   Assume that
 $X$ has isolated  singularities,
 $-K_X$ is $\q$-Cartier and $f$-ample, and  $Z$ is reduced and
geometrically irreducible.
 We would like to explain what it means that $Z$ is a (possibly
nonnormal) Del Pezzo  surface.

Let $\Sigma$ be the finite set of  points of $X$  where $K_X$ is not
Cartier. Set
$X^0:=X\setminus \Sigma$ and 
$Z^0:=Z\setminus \Sigma$.   (To be completely precise, I should let
$\Sigma$ be the set of points where $X$ is not Gorenstein. At the end,
however, this does not make any difference.)
$Z^0$ is a Cartier divisor on the   variety $X^0$, thus 
$$
\omega_{Z^0}\cong \omega_{X^0}(Z^0)\otimes \o_{Z^0}\cong 
\omega_{X^0}\otimes \o_{Z^0}.
$$ By assumption  $\omega_{X^0}^{-m}$ is $f$-very ample for $m\gg 1$,
so we obtain that $\omega_{Z^0}^{-m}$ is  very ample for $m\gg 1$. 

Let $p:\bar Z\to Z$ be the normalization and set $\bar
Z^0:=\pi^{-1}(Z^0)$. There is an adjunction map
$\omega_{\bar Z^0}\to \pi^*\omega_{Z^0}$ (cf.\
\cite[Ex.III.7.2.a]{Hartshorne77}). (The adjunction map is defined even
over $\bar Z$ but $\pi^*\omega_Z$ may be messy since
$\omega_Z$ need not be locally free.) $\bar Z^0$ is a normal surface,
so there is an effective  divisor $B^0\subset \bar Z^0$ such that
$$
\o_{\bar Z^0}(K_{\bar Z^0}+B^0)\cong \omega_{\bar Z^0}(B^0)\cong
\pi^*\omega_{Z^0}.
$$  (Further information about $B^0$ is given in (\ref{what-is-B?}).)
Here both sides are locally free (since the right  hand side is), so we
can raise them to any power to get isomorphisms
$$
\o_{\bar Z^0}(m(K_{\bar Z^0}+B^0))\cong 
\pi^*\omega^m_{Z^0}.
$$ Now choose $m>0$ such that $mK_X$ is Cartier. The reflexive  sheaves
$$
\o_{\bar Z}(-m(K_{\bar Z}+B))
\qtq{and}
\pi^*\o_X(-mK_X)
$$ are isomorphic on $\bar Z^0$, thus outside finitely many points.
Hence they are isomorphic. Since $\o_X(-mK_X)$ is locally free we
conclude that
\begin{enumerate}
\item $K_{\bar Z}+B$ is $\q$-Cartier and $-(K_{\bar Z}+B)$ is ample, and
\item $(K_{\bar Z}+B)^2=K_{X_t}^2$ where $X_t$ is any smooth fiber of
$X/C$.
\end{enumerate} We can also get some information about the global
sections. One has to be a little careful since $K_X$ may not be Cartier
everywhere. As usual, let $\omega_X^{[r]}$ denote the double dual of
$\omega_X^{\otimes r}$.
$\omega_X^{[-n]}$ is reflexive, so 
$\omega_X^{[-n]}\otimes \o_Z$ has no embedded points. Thus
\begin{eqnarray*} H^0(Z^0,\omega^{-n}_{Z^0})&=&
H^0(Z^0,\omega_X^{[-n]}\otimes \o_{Z^0})\supset
H^0(Z,\omega_X^{[-n]}\otimes \o_Z)\\ &\geq &
H^0(X_t,\omega_{X_t}^{[-n]})={\textstyle\binom{n+1}{2}} (K_{X_t}^2)+1,
\end{eqnarray*} where the inequlity holds by semicontinuity of $H^0$
and the last equality is a fact about smooth Del Pezzo surfaces. We can
pull back these sections to $\bar Z^0$ and  extend them across the
finitely many points. Thus we conclude that
\begin{enumerate}\setcounter{enumi}{2}
\item $H^0({\bar Z},\o_{\bar Z}(-(K_{\bar Z}+B)))\geq (K_{X_t}^2)+1\geq
2$.
\end{enumerate} Choose a pencil in  $|-(K_{\bar Z}+B)|$ and write it as
$D+M$ where $D$ is the fixed part and $M$ the moving part.  If   $K_X$
is Cartier at all points of $X(\r)$ then $Z(\r)\subset Z^0$ and also
$\bar Z(\r)\subset \bar Z^0$. This shows that $K_{\bar Z}+B$ is Cartier
at  all real points. Setting  $S:=\bar Z$  we obtain a quadruplet
$(S,B,D,M)$ and we have proved the following.
\end{say}

\begin{prop} Let $f:X\to C$ be a family of Del Pezzo surfaces
satisfying the conditions (\ref{bas.setup.not}.1--3). Let $Z:=f^{-1}(0)$
be a reduced and geometrically irreducible fiber. Then the  quadruplet
$(S,B,D,M)$ constructed above satisfies the conditions
(\ref{S-conds}).\qed
\end{prop}

\begin{condition}\label{S-conds}
  For  a quadruplet 
$(S,B,D,M)$ consider the following properties.
\begin{enumerate}
\item $S$ is a normal projective surface over $\r$,
\item  $B$ and $D$ are  Weil divisors on $S$,
\item $M$ is a pencil on $S$ without fixed components,
\item $K_S+B+D+M\sim 0$.
\item $K_S+B$ is Cartier at all real points of $S$,  and
\item $-(K_S+B)$ is $\q$-Cartier and ample.
\end{enumerate}
\end{condition}

\begin{rem}  (1) Most quadruplets satisfying the above conditions do
not arise as special fibers of  families of Del Pezzo surfaces. An
obvious numerical conditions is that $(K_S+B)^2$ be an integer beween 1
and 9.  Even if this holds, there is no reason to assume that the
surface  $Z$ is smoothable. 

(2) Even if we know $(S,B,D,M)$, it is not always easy to reconstruct
the fiber $Z$. First of all, we need to   know the map $B\to
(\mbox{singular curve of $Z$})$;  a set theoretic information. If $B$
is reduced then this determines the scheme structure of $Z$ at least
when $Z$ is $S_2$ (which holds if $X$ has terminal singularities). If,
however, $B$ is not reduced,  additional scheme theoretic
information is needed and this seems rather complicated. See
\cite{Reid94} for a closely related case.
\end{rem}

\begin{say}\label{what-is-B?} Let $X$ be a smooth 3--fold and $Z\subset
X$ a  reduced surface which is not normal along a curve $C\subset Z$.
Let $\pi:\bar Z\to Z$ be the normalization and set
$\red(\pi^{-1}(C))=\sum_i\bar C_i$. As in (\ref{degen.of.dp.say}), we
see that
$\omega_{\bar Z}(\sum b_i\bar C_i)\cong \pi^*\omega_Z$ for some
$b_i\geq 0$. We would like to establish a relationship between the
coefficients $b_i$ and the local structure of $Z$ along $C$.

We can cut everything by a general hyperplane, and our problem is
reduced to a curve question which has been classically studied in
detail, since
$(\sum b_i)/2$ is exactly the contribution of the singular point to the
arithmetic genus of a curve. We need the classification of
singularities with small $b_i$:
\begin{enumerate}
\item $b_1=0$ iff $Z$ is smooth along $C$.
\item $b_1=1$ iff $Z$ has 2 branches meeting transversely along $C$.
\item $b_i\leq 2$ for every $i$ iff either  $Z$ has one branch with a
cusp along
$C$, or 2 branches which are simply tangent along $C$ or  3 branches  
meeting pairwise transversely along $C$ or 
$b_i\leq 1$ for every $i$.
\end{enumerate}
\end{say}

\begin{say}\label{S^m-defn.say}
 In general $S$ is singular and so we try to study it through a
suitable resolution. The most natural choice would be to take its
minimal resolution, but the following partial resolution turns out to
be  more convenient.  (The main reason for this choice is explained  in
(\ref{top.go.down.rem}.3).) 

 There is a unique morphism $f:S^m\to S$  such that 
\begin{enumerate}
\item $f$ is  an isomorphism above $P\in S$ if $S$ is smooth at $P$ or
if $S$ has a Du Val singularity at $P$ and $P\not\in \supp B$
\item  $f$ is  the minimal resolution over   $P\in S$ otherwise.
\end{enumerate} (Note that  this is not the same as the so called
``minimal Du Val resolution".)

$f$ is a birational map between normal surfaces, thus we can pull back
any divisor by $f$ if we allow rational coefficients. Hence we can
define $B^m$ by the formula $K_{S^m}+B^m\equiv f^* (K_S+B)$. 
$-B^m\equiv_f K_{S^m}$ has nonnegative intersection number with any
exceptional curve, hence $B^m$ is effective (cf.\
\cite[3.39]{KoMo98}). Moreover, every exceptional curve appears in
$B^m$ with positive coefficient  (cf.\ \cite[4.3--5]{KoMo98}). (Here we
use that we did not resolve Du Val points.) Write $f^*M=D'+M^m$ where
$D'$ is the fixed part (which may have rational coefficients) and $M^m$
is a pencil without fixed components.
 Set $D^m:=f^*D+D'$. 

\begin{enumerate}
\item[$(*)$] In this paper, the quadruplet
$(S^m,B^m,D^m,M^m)$   always denotes the one constructed above  staring
with  $(S,B,D,M)$.
\end{enumerate}
\end{say}

\begin{condition}\label{S^m-conds}
 For a quadruplet
$(S^m,B^m,D^m,M^m)$  consider the following properties.
\begin{enumerate}
\item $S^m$ is a  projective surface over $\r$ with Du Val
singularities,
\item  $B^m$ and $D^m$ are effective $\q$-divisors on $S$ such that
$B^m+D^m$ is an integral divisor,
\item $M^m$ is a pencil on $S^m$ without fixed components,
\item $K_{S^m}+B^m+D^m+M^m\sim 0$.
\item If $C\subset S^m$ is a geometrically irreducible real curve then
$C$ appears in $B^m$ with integer coefficient,
\item $-(K_{S^m}+B^m)$ is nef and big
\item $-(K_{S^m}+B^m)$  has positive intersection number with every
$(-1)$-curve (\ref{-1.curve.def}),
\item $S$ is smooth along $\supp B$,
\item $-(K_{S^m}+B^m)$  is nef and has positive intersection number
with every curve not in $\supp B$.
\end{enumerate}
\end{condition}

\begin{prop}  If $(S,B,D,M)$ satisfies the conditions (\ref{S-conds})
then  $(S^m,B^m,D^m,M^m)$ satisfies the conditions (\ref{S^m-conds}).
\end{prop}

Proof. (1,3,4,6,8) are clear from the construction.
$B^m+D^m\sim -K_{S^m}-M^m$ now implies (2). If 
$C\subset S$ is a geometrically irreducible real curve and
$f(C)$ is a curve then  $C$ appears in $B$ and in $B^m$ with the same
integer coefficient. If $f(C)=P$ is a point then $P$ is real, so
$K_S+B$ is Cartier  there and every curve above $P$ appears with
integral coefficient in $B^m$. 

(9) holds since every $f$-exceptional  curve is in $B$. Since we took
minimal resolutions, there are no $f$-exceptional  $(-1)$-curves, and
this shows (7).\qed
\medskip

\begin{lem}\label{1-8.conds.descend}
 Let $(F,B,D,M)$ be a quadruplet.  Let $g:F\to F'$ be a birational
morphism and set $B':=g_*B,D':=g_*D$ and $M':=g_*M$. Then:
\begin{enumerate}
\item  If  $(F,B,D,M)$  satisfies one of the conditions
(\ref{S^m-conds}.2-6 or 9)  then  $(F',B',D',M')$ also satisfies the
same condition.
\item If $g^{-1}$ is a composite of $(1,m)$-blow ups and $(F,B,D,M)$ 
satisfies one of the conditions (\ref{S^m-conds}.1-9)  then 
$(F',B',D',M')$ also satisfies the same condition.
\end{enumerate}
\end{lem}

Proof.  This is clear for (\ref{S^m-conds}.2,3 and 5) and
$K_{F'}=g_*K_F$ implies it for  (\ref{S^m-conds}.4). 

Let $E_i$ be the connected components of the exceptional set of $g$. 
We can write $K_F+B\equiv g^*(K_{F'}+B')+E$ where $E$ is
$g$-exceptional. $-(K_F+B)$ is nef, so  $E$ is effective and
 $E_i\subset \supp E$  unless $K_F+B$ is numerically trivial along
$E_i$ (cf.\ \cite[3.39]{KoMo98}).  If $C\subset F$ is an  irreducible
curve  which is not
$g$-exceptional then
$$ g_*(C)\cdot (K_{F'}+B')= C\cdot g^*(K_{F'}+B')= C\cdot  (K_F+B)
-C\cdot E\leq C\cdot  (K_F+B).
$$ This shows that $-(K_{F'}+B')$ is nef if $-(K_{F}+B)$ is, settling
the case  (\ref{S^m-conds}.9).  If (\ref{S^m-conds}.6) holds on $F$ then
bigness of $-(K_{F'}+B')$  follows from
$$ (K_{F'}+B')^2=(K_F+B)\cdot g^*(K_{F'}+B') =(K_F+B)^2-(K_F+B)\cdot
E\geq (K_F+B)^2.
$$

Assume now that $g^{-1}$ is a single $(1,m)$-blow up.  The points
$g(E_i)$ are smooth on $F'$  by (\ref{(1,m)-bu.defn}), which implies
the claim for  (\ref{S^m-conds}.8). Every connected component of the 
exceptional set is a
$(-1)$-curve by (\ref{(1,m)-bu.defn}) and $K_F+B$ has negative
intersection number with it. So $\supp E$ coincides with the
exceptional set and  if $C$ intersects the exceptional set then 
$$ g(C)\cdot (K_{F'}+B')\leq -C\cdot E<0.
$$ This shows the claim for (\ref{S^m-conds}.7). By induction this gives
(\ref{1-8.conds.descend}.2).\qed

\section{The Classification of $S^r$}

\begin{notation}\label{S^r.class.not}
 Throughout this section 
$(S^m,B^m,D^m,M^m)$ denotes a  quadruplet which satisfies the
conditions  (\ref{S^m-conds}.1--8).  Let $g:S^m\to S^*$  be a minimal
model of $S^m$ and set $B^*:=g_*(B^m)$.  For simplicity we assume $S^*$
is not obtained by $(1,m)$-blow ups from another surface. (This could
happen in only a handful of cases. For instance, the blow up of $\p^2$
at  a point is also a 
$\p^1$-bundle so it could be $S^*$ as in (\ref{surf.dv.mmp}.4.C).)
\end{notation}

\begin{say}[How to determine $S^m$?]\label{backwards.meth}{\ }

We use the following method to get information about $S^m$.

(1) First we determine the possible surfaces $S^*$. This was in fact
done using the MMP for real surfaces with Du Val singularities. 

(2) Then we get a list of all possible quadruplets $(S^*,B^*,D^*,M^*)$.
This involves finding all possible ways of writing
$-K_{S^*}\sim (B^*+D^*)+M^*$, so this is equivalent to classifying all
pencils in $|-K_{S^*}|$.  This is easy to do if $|-K_{S^*}|$ is small.
Unfortunately, $|-K_{S^*}|$ gets arbitrarily large  for minimal
ruled surfaces  and I do not know of any useful classification of all
anticanonical pencils. 

(3) Given a quadruplet $(S^*,B^*,D^*,M^*)$ we can try to find all
possible ways it came from an $(S^m,B^m,D^m,M^m)$. We factor
$S^m\to S^*$ into $(1,m)$-blow ups
$$ S^m\to \cdots \to S_2\to S_1\to S^*.
$$ For any of the intermediate stages $g:S^m\to S_i$ set
$B_i:=g_*B^m,D_i:=g_*D^m,M_i:=g_*M^m$.   A key point to observe is that
by (\ref{1-8.conds.descend}) all the  
$(S_i,B_i,D_i,M_i)$ satisfy the conditions (\ref{S^m-conds}.1--8). Thus
we can work our way backwards one blow up at a time starting with
$S^*$.

(4) Assume that we already have $(S_i,B_i,D_i,M_i)$ and that
$\pi:S_{i+1}\to S_i$ is the ordinary blow up of a point $P\in S_i$ with
exceptional curve $E\subset S_{i+1}$. By (\ref{S^m-conds}.4) 
$$ K_{S_{i+1}}+B_{i+1}+D_{i+1}+M_{i+1}=
\pi^*(K_{S_{i}}+B_{i}+D_{i}+M_{i}),
$$ so  we conclude that the
$$
\mbox{coefficient of $E$ in $(B_{i+1}+D_{i+1})$}=
\mult_P(B_{i}+D_{i})+\mult_PM_i-1.
$$
  In particular, we can blow up only points in $\supp(B_{i}+D_{i})$ and
the base points of $M_i$.

(5) Given $(S^*,B^*,D^*,M^*)$ and the sequence of blow ups leading to
$S^m$, we have determined $S^m,B^m+D^m$ and $M^m$.  Conditions 
(\ref{S^m-conds}.5--8) give further restrictions on $B^m$. In many
cases these are impossible to satisfy.  

(6) The role of (\ref{S^m-conds}.5) turns out to be crucial for us. 
Frequently there are many possibilities for $B^m$ such that 
$(S^m,B^m,D^m,M^m)$ satisfies all the conditions (\ref{S^m-conds})
except (\ref{S^m-conds}.5), but there are  no choices of $B^m$ where
(\ref{S^m-conds}.5) also holds.

(\ref{S^m-conds}.5) is especially useful if there are many 
geometrically irreducible  curves in $\supp(B^m+D^m)$.  This is the
reason for studying $S^r$ since  in $S^r\to S^*$ all exceptional curves
are geometrically irreducible.  It   turns out that $S^r$ is obtained
from $S^*$ by at most 1 blow up and this allows us to understand  all
possible
$(S^r,B^r,D^r,M^r)$  quite well.

(7)  In (\ref{5-cases.lem})  we subdivide the possible pairs
$(S^*,B^*)$ into 5 cases and then do a separate classification in each
case.
\end{say}

\begin{lem}\label{5-cases.lem}
 With the above notation the pair $(S^*,B^*)$ satisfies one of the
following conditions.
\begin{enumerate}
\item  $S^*$ is a Del Pezzo surface with $\rho(S^*)=1$ which is neither
$\p^2$ nor a quadric in $\p^3$.
\item  $S^*$ is a weak Del Pezzo surface with $\rho(S^*)=2$  and
$B^*=0$.
\item    There is a $\p^1$- bundle  structure $S^*\to A$ such that 
 $B^*_{\c}$ has a unique irreducible component   which dominates $A$.
This component is a section.
\item   There is a conic bundle structure $S^*\to A$ such that 
$B^*_{\c}$ has two irreducible components   which dominate $A$. These
components are  conjugate  sections.
\item  $S^*$ is   $\p^2$ or a quadric in $\p^3$.
\end{enumerate}
\end{lem}

Proof.  If $S^*$ is a Del Pezzo surface with $\rho(S^*)=1$ then we have
either  (1) or (5). In all other cases there is a conic bundle
structure $S^*\to A$. Let $F\subset S^*$ be a general fiber. Then 
$$ 2=-K_{S^*}\cdot F=(B^*+D^*)\cdot F+M^*\cdot F\geq (B^*+D^*)\cdot F.
$$
$B^*+D^*$ is a sum curves with integral coefficients,  so it intersects
$F$ in at most 2 points. Thus $B^*_{\c}$ has at most  2 horizontal
components.

$-(K_{S^*}+B^*)\cdot F>0$ since  $-(K_{S^*}+B^*)$ is big, thus
$B^*\cdot F<2$.

If $B^*_{\c}$ has 1 horizontal component  $H$, then $H$ is defined over
$\r$ and  hence $H$ appears in $B^*$ with integer coefficient.  Together
with
$B^*\cdot F<2$  this shows that  $H$ is a section. Then $S^*$ is
 a $\p^1$-bundle  by (\ref{cb.with.sect}) and (\ref{S^m-conds}.8), so we
are in  case (3).

If  $B^*_{\c}$ has 2 horizontal components and they are both defined
over $\r$ then arguing  as above leads to a contradiction.

If $B^*_{\c}$ has 2 horizontal components which are conjugates  $H,
\bar H$ then they are both sections and we are in case (4). 

If    $B^*=0$ then we are in  case (2).

Finally it may happen that $B^*_{\c}\neq 0$ has no horizontal
components.  $\rho(S^*)=2$ and so $\rho(S^*/A)=1$. $S^*$ is
$\q$-factorial  hence every curve in $S^*$ is $\q$-Cartier. This
implies that 
 any vertical curve is numerically equaivalent to a (rational) multiple
of a general fiber, in particular $B^*$ is nef. This implies that
 $-K_{S^*}$ is ample. By the Cone Theorem  (cf.\ \cite[3.7]{KoMo98})
 $S^*$ has another extremal ray with contraction
$S^*\to W$. This  is not a birational contraction  by assumption
(\ref{S^r.class.not}). So   $S^*\to W$ is another conic bundle
structure and  the vertical $B^*_{\c}$ becomes horizontal for $S^*\to
W$. These cases were treated already.
\qed

\begin{thm}\label{dp-case.thm}
 Assume that $S^*$ is a weak Del Pezzo surface  with Du Val
singularities, not isomorphic to $\p^2$ or  to a quadric   in
$\p^3$. Assume furthermore that either $\rho(S^*)=1$ or $B^*=0$.  Then 
 $S^r$ is a weak Del Pezzo surface  with Du Val singularities and
$B^r=0$. 
\end{thm}

Proof. Assume first  that
$B^*\neq 0$. Write $B^*+D^*=C+C'$ where $\supp C=\supp B^*$. 
$-K\equiv C+C'+M^*$ and  $C$ is a Cartier divisor by (\ref{S^m-conds}.8)
which is not empty if $B^*\neq 0$.  On $S^*$ every effective Cartier
divisor  is ample since $\rho(S^*)=1$, so
$C'+M^*$ is ample. This is impossible by (\ref{dp.pic.divis}). 

Thus  $B^*= 0$ and $S^*\not\cong \p^2$. If $B^r= 0$ then
$-(K_{S^r}+B^r)=-K_{S^r}$ is nef and big,
 so $S^r$ is a  weak Del Pezzo surface  with Du Val singularities.

We need to exclude the case when $B^r\neq 0$. Let 
$$
\cdots \to S_2\to S_1\to S^*
$$
 be the series of blow ups leading to
$S^r$ and  $S_{i+1}\to S_i$  the last blow up  with exceptional curve
$E\subset S_{i+1}$ such that 
$B_i=0$. Then $B_{i+1}=aE$ for some $a>0$ and $a$ is an integer since
$E$ is geometrically irreducible. $S_{i+1}$ is smooth along
$B_{i+1}$, so $S_{i+1}\to S_i$ is an ordinary blow up. This is  
impossible by (\ref{0<-B.imposs}).\qed

\begin{prop}\label{0<-B.imposs}
 Let $F$ be a projective surface with Du Val singularities over
$\c$ and  
$p:G\to F$   the blow up of a smooth point with exceptional curve
$E\subset G$. Assume that $ -(K_G+aE)$ is nef and has positive
intersection number with every $(-1)$-curve for some $a\geq 1$.

Then $F\cong \p^2$ and $a<2$.
\end{prop}

Proof. Let $F'\to F$ and $G'\to G$ be the minimal resolutions. Then
$F'$ is the blow up of $G'$ at a point and $ -(K_{G'}+aE')$ is the pull
back of $ -(K_F+aE)$. Thus it is sufficient to consider the case when
$F$ itself is smooth.

If $F\not\cong\p^2$ then there is a morphism $g:F\to \p^1$ whose
general fiber is $\p^1$. $-K_F$ is $g$-nef, so it is easy to see that
every fiber of
$g$ has dual graph
$$
\stackrel{0}{\circ} \qtq{or} 
\stackrel{-1}{\circ} - 
\stackrel{-2}{\circ} - \stackrel{-2}{\circ} - 
 \cdots - 
\stackrel{-2}{\circ} - \stackrel{-1}{\circ}.
$$ If $P$ is on a fiber of the second type then 
$E$ intersects a rational curve with self intersection $-2$ or $-3$ and
$ -(K_G+aE)$ is not nef for $a>0$. In the first case, the birational
transform of the fiber becomes a $(-1)$-curve $F$ 
 intersecting $E$. Thus $-(K_G+aE)\cdot F=1-a\leq 0$.\qed

\begin{thm}\label{p1-bund-case.thm} Assume that $S^*\to A$ is a
$\p^1$-bundle and 
$H\subset B^*$ is a (real) section. Then 
 $S^r=S^*$. 
\end{thm}

Proof.  If $S^r\neq S^*$ then there is a last contraction
$S_1\to S^*$. The inverse of this is a $(1,m)$-blow up of a real point
 $P\in S^*$ and we need to show that this can not happen. Let $F'$ be
the fiber of $S^*\to A$ containing $P$. 

If $S_1\to S^*$ is an ordinary blow up then the exceptional curve
$E\subset S_1$ and the birational transform $F$ of $F'$  are both
$(-1)$-curves intersecting in a point. By (\ref{S^m-conds}.7) this
implies that
$E\cdot(K_{S_1}+B_1)<0$ and
$F\cdot(K_{S_1}+B_1)<0$. We can write $B_1=b_eE+b_fF+B'$ where $B'$
does not have $E, F$ as components. $H$ intersects either $E$ or $F$, so
$E\cdot B'\geq 1$ or  $F\cdot B'\geq 1$. We get a contradiction by
(\ref{ineq.lem.1}). 

The $(1,m)$-blow ups are excluded by  (\ref{(1,m)-bu.say}).\qed

\begin{say}[Excluding $(1,m)$ blow ups]\label{(1,m)-bu.say}

Given a pair $(S^*,B^*)$ assume that  $S_1\to S^*$ is a $(1,m)$-blow up
for some $m\geq 2$ with exceptional curve $E$. Then 
$S_1$ is singular at a point   $P\in E$, so $E$ is not in $B_1$. Let
$\pi:\tilde S_1\to S_1$ be the  minimal resolution of $P$. Write
$\pi^*(D_1+M_1)=\tilde D_1+\tilde M_1$ where $\tilde M_1$ has no fixed
components.  Then $(\tilde S_1,\pi^*B_1,  \tilde D_1,\tilde M_1)$
satisfies the conditions (\ref{S^m-conds}.1--8) and $\tilde S_1$ is
obtained from $S^*$ by $m$ ordinary blow ups (\ref{(1,m)-bu.defn}). So
once we prove that there can be at most one ordinary blow up, this
implies that there are no
$(1,m)$-blow ups at all for $m\geq 2$. 
\end{say}

\begin{lem}\label{ineq.lem.1}
 Let  $S$ be a surface, $E,F\subset S$ effective curves and $B'$ an
effective $\q$-divisor on $S$ whose support does not contain $E$ and
$F$.  Assume that $(E^2)=(F^2)=(E\cdot K_S)=(F\cdot K_S)=-1$ and
$(E\cdot F)=1$. Let $b_e,b_f$ be integers and assume that
$$ E\cdot (K_S+b_eE+b_fF+B')<0\qtq{and} F\cdot (K_S+b_eE+b_fF+B')<0.
$$ Then $b_e=b_f$, $E\cdot B'<1$ and $F\cdot B'<1$.\qed
\end{lem}

\begin{thm}\label{cb-bund-case.thm}
 Assume that $S^*\to A$ is a conic bundle and  $B^*_{\c}$ contains a
pair of  conjugate  sections. Then  $A$ is rational and  either 
$S^r=S^*$ or $S^r$ is obtained from $S^*$ by one ordinary blow up.
\end{thm}

Proof.  Let $H,\bar H\subset \supp B^*$ be the conjugate sections.
First we determine $B^*+D^*$. We can write $-K_{S^*}\equiv H+\bar H+aF$
for some $a$ where $F$ is a general fiber. $B^*+D^*$ is a Weil divisor
containing $H+\bar H$, so
$(B^*+D^*-H-\bar H)+M\equiv aF$ is effective and moves in a pencil,
hence
$a\geq 1$. By the adjunction formula
$$ 2g(H)-2=H\cdot (K_{S^*}+H)=-(H\cdot \bar H)-a.
$$ Thus $g(H)=0$ and either $H\cdot \bar H=1$ and  $B^*+D^*=H+\bar H$
or  $H\cdot \bar H=0$ and  $B^*+D^*=H+\bar H+\epsilon C$ where
$\epsilon\in \{0,1\}$ and $C$ is a fiber.  In the first case let $C$
denote the fiber passing through the unique point of  $H\cap \bar H$.
In both cases $M=|F|$ is base point free. The   real points of $H+\bar
H+\epsilon C$ are in  $C(\r)$. So all real blow ups take place over $C$
by (\ref{backwards.meth}.4). We distinguish three cases according to
the type of
$C$.

(1): $C$ is a double fiber. Any section intersects a double fiber at a
singular point of $S^*$. Since $H\subset B^*$ and $B^*$ is disjoint from
the singular points of $S^*$, we see that $S^*\to A$ does not have any
double fibers.

(2): $C$ is a smooth fiber. Consider   the case of 2 ordinary real blow
ups $S_2\to S_1\to S^*$. There are 3 cases and after blow up we get one
of the following  curve configurations where $\bullet$ denotes the
birational transform of $C$.
$$
\stackrel{-1}{\circ} - \stackrel{-2}{\bullet} - 
\stackrel{-1}{\circ}
\qtq{or}
\stackrel{-1}{\circ} - \stackrel{-2}{\circ} - 
\stackrel{-1}{\bullet}
\qtq{or}
\stackrel{-2}{\circ} - \stackrel{-1}{\circ} - 
\stackrel{-2}{\bullet}.
$$ Let $E_1,E_2,E_3$ be these curves from left to right.
$B_2=e_1E_1+e_2E_2+e_3E_3+H_2+\bar H_2$, the
$e_i$ are integers, $-(K_{S_2}+B_2)$ is positive on the $(-1)$-curves
and nonnegative on the $(-2)$-curves. By solving these inequalities we
obtain that
$H_2+\bar H_2$ does not intersect the $(-2)$-curves,
$e_1=e_2=e_3$ in the first 2 cases  and
 $e_1=e_2/2=e_3$ in the third case.

If $C\cap (H+\bar H)$ is a conjugate point pair, then 
$H_2+\bar H_2$   intersects the curve $\bullet$ and this leads to a
contradiction in the first and third cases.  In the second case  we use
(\ref{backwards.meth}.4) to conclude that 
$B^*+D^*=H+\bar H+C$ and   $e_2=0$,  giving a contradiction.

If $C\cap (H+\bar H)$ is a real point $P$, then $C$ has coefficient
zero in $B^*+D^*$ hence we get that $e_1=e_2=e_3=0$ and  the only point
that we can blow up on $S^*$ is $P$. Thus the first case can not happen
and $H_2+\bar H_2$ intersects $E_2$ in the second case and 
$E_1$ in the third  case.   Both are impossible.

 $(1,m)$-blows up are   excluded by  (\ref{(1,m)-bu.say}).

(3): $C$ is a reducible fiber.  The only real point of $C$ is its
singular point, hence we have to start by blowing it up. There is only
1 case and after 2 blow ups we get   the following  curve configuration.
$$
\stackrel{-1}{\circ} - \stackrel{-2}{\circ} -  {\bullet},
$$ where $\bullet$ denotes the birational transform of $C$ (thus it is
geometrically reducible). Let $E_1,E_2,E_3$ be these curves from left
to right.
$B_2=e_1E_1+e_2E_2+e_3E_3+H_2+\bar H_2$, $e_1,e_2$ are integers,
$-(K_{S_2}+B_2)$ is positive on the $(-1)$-curve and nonnegative on the
$(-2)$-curves. By solving these inequalities we obtain that
$e_1=e_2=e_3$ and 
$H_2+\bar H_2$ does not intersect   $\bullet$. This is  impossible
since $H$,  a section, does not pass through the singular point of the
fiber $C$.

$(1,m)$-blow ups are excluded as above.
\qed

\begin{notation} There are 5 normal quadrics over $\r$ up to
isomorphism. The smooth ones  are $Q^{4,0}:=(x^2+y^2+z^2+t^2=0)$,
$Q^{3,1}:=(x^2+y^2+z^2-t^2=0)$ and
$Q^{2,2}:=(x^2+y^2-z^2-t^2=0)$. The quadric cones are
$Q^{3,0}:=(x^2+y^2+z^2=0)$ and $Q^{2,1}:=(x^2+y^2-z^2=0)$.

$Q^{4,0}(\r)=\emptyset$ and $Q^{3,0}(\r)$ is a single point, so they
are not very interesting for us.
\end{notation}

\begin{thm}\label{p2-case.thm}
 Assume that $S^*$ is  $\p^2$ or a quadric in $\p^3$. Then  one of the
following holds.
\begin{enumerate}
\item   $S^r$ is  $\p^2$ or a quadric in $\p^3$. 
\item $S^r$ is weak Del Pezzo surface and $B^r=0$. 
\item $S^r$ is obtained from $\p^2$ by one ordinary blow up.
\item $S^r$ is $\p^2$ blown up in 2 points and $B^r(\r)=\emptyset$.
\item $S^r$ is one of the above cases blown up at conjugate pairs of
points.
\end{enumerate}
\end{thm}

(The last case occurs since blowing up the quadric
$Q^{3,1}$ at a real point is the same as blowing up
$\p^2$ at a conjugate pair of points. A more economical choice of $S^*$
 would have eliminated this case.)
\medskip

Proof.  Assume that we blow up 2 (possibly infinitely near) real points
in
$\p^2$. There is a unique line passing through the center of both blow
ups. This line and the two exceptional curves form a configuration
$$
\stackrel{-1}{\circ} - \stackrel{-1}{\circ} - 
\stackrel{-1}{\circ}
\qtq{or}
\stackrel{-1}{\circ} - \stackrel{-1}{\circ} - 
\stackrel{-2}{\circ}.
$$ Let $E_1,E_2,E_3$ be these curves from left to right.
$B_2=e_1E_1+e_2E_2+e_3E_3+B'_2$, the
$e_i$ are integers, $-(K_{S_2}+B_2)$ is positive on the $(-1)$-curves
and nonnegative on the $(-2)$-curve. By solving these inequalities we
obtain that
$e_1=e_2=e_3=0$ and $E_i\cdot B'_2<1$ for every $i$. 

If $B_2=0$ then we end up in case (2) as in (\ref{dp-case.thm}).  If
$B^*$ contains a   curve with positive integer coefficient then its
birational transform intersects $E_1+E_2+E_3$ and one of the
inequalities $E_i\cdot B'_2<1$ is violated.

A similar argument excludes the possibility of even one blow up when
$S^*$ is  a quadric  and $B^*$ contains a   curve with positive integer
coefficient.

We are left with three cases: 
$S^*=\p^2$, $B^*=c(L+\bar L)$ for a conjugate pair of lines, or
$S^*=Q^{3,1}$, $B^*=c(L+\bar L)$ for a conjugate pair of intersecting
lines or 
$S^*=Q^{2,2}$, $B^*=c(L+\bar L)$ for a conjugate pair of
nonintersecting lines. In the latter case $S^*\cong \p^1\times \p^1$
and one of the projections lands us in case (\ref{cb-bund-case.thm}), so
this is already treated.

Assume  that $S^*=\p^2$  and $B^*=c(L+\bar L)$ for a  conjugate pair of
lines intersecting at a (necessarily real) point $P$. Then
$B^*+D^*=L+\bar L$ and
$M$ is a pencil of lines with a base point $Q$. By
(\ref{backwards.meth}.4) the only possibilities to blow up are $P$ and
$Q$. If $Q\neq P$ and we blow up
$Q$ then projection from $Q$ becomes a $\p^1$-bundle and $L,\bar L$
become conjugate sections. This case was treated in 
(\ref{cb-bund-case.thm}) and we get that blowing up $P$ is the only
possible further blow up.

If we blow up $P$ (possibly $P=Q$) and the exceptional curve  $E$
appears in
$B$ with coefficient $\geq 1$ then we are in case
(\ref{p1-bund-case.thm}) and no more blow ups are possible. So $E$
appears in $B$ with coefficient 0. Another blow up on $E$ would create
a $(-2)$-curve with coefficient 0 in $B_2$ intersecting the birational
transform of $L$, a contradiction. Thus we can blow up only $Q\neq P$
and then we are in the already discussed case. 

As we showed, we can never blow up infinitely near points, so there are
no $(1,m)$-blow ups by (\ref{(1,m)-bu.say}).

If $S^*=Q^{3,1}$  then a one point blow up of $S^*$ is also a   blow up
of $\p^2$ at a  conjugate pair of points. So we are  reduced to
considering the blow ups of
$\p^2$. These cases can also  be treated directly.
\qed
\medskip

We can summarize the above results as follows.

\begin{prop}\label{S^r-summ.prop}
 Assume that the quadruplet $(S^m,B^m,D^m,M^m)$ satisfies the
conditions (\ref{S^m-conds}). Then $(S^r,B^r)$ satisfies one of the
following conditions.
\begin{enumerate}
\item $S^r$ is a weak Del Pezzo surface and $B^r=0$.
\item  $S^r$ is a $\p^1$-bundle, $B^r$ consists of a section, at most 2
real fibers and some conjugate pairs of fibers.
\item $S^r$ is a  conic bundle, $B^r$ consists of a conjugate pair of
sections and at most one fiber.
\item $S^r$ is a  conic bundle with one point blown up, $B^r$ consists
of a conjugate pair of sections and possibly one fiber which is the
union of a $\p^1$ and of a conjugate pair of rational curves.
\item $S^r=\p^2$ and $B^r$ is  either a line or a pair of lines or a
conic.
\item $S^r=Q^{2,2}\cong \p^1\times \p^1$ and $B^r$ is   either   a line
or a pair of intersecting lines or a conic.
\item $S^r=Q^{2,1}$ and  $B^r$ is a hyperplane section not through the
vertex.
\item  $S^r(\r)$ is a finite set.
\end{enumerate}
\end{prop}

Proof. There are only a few points which have not been settled earlier.

In the $\p^1$-bundle case write $B^r=H+B'$. From the adjunction formula
we get that
$$ -2\leq 2g(H)-2=H\cdot (K+H)=-H\cdot (B'+D^r+M^r)\leq -(H\cdot B'),
$$ so $B^r$ contains at most 2 fibers.

In the conic bundle case we either have no blow ups and get (3) or let
$C$ denote the fiber containing the center of the blow up. There are 2
cases to consider. $C$ can be singular giving case (4). If $C$ is
smooth then after one blow up we get 2 intersecting
$(-1)$-curves and both  appear with coefficient 0 in $B^r$.\qed
\medskip

The following result  could have been proved much earlier.

\begin{lem} Let 
$(S,B,D,M)$  be a quadruplet 
 satisfying the conditions (\ref{S-conds}). Then  one of the following
holds.
\begin{enumerate}
\item $S$ is a cone over an elliptic curve and
$B=0$.
\item  $S_{\c}$ is a rational surface with rational singularities  and
every irreducible component of $B^m$ is a smooth rational curve.
\end{enumerate}
\end{lem}

Proof. In the previous classification, a  nonrational surface can occur
only in (\ref{p1-bund-case.thm}). So $S^*$ is a $\p^1$-bundle and
$H\subset B^*$ is a section. From the adjunction formula
$$ 2g(H)-2=H(K+H)=-H(B^*-H+D^*+M^*)\leq 0.
$$ So $H$ is elliptic, $B^*=H$ and $H(K+B^*)=0$. We have proved that
there are no real blow ups and a similar computation shows that there
are no complex blow ups either. So $S^m=S^*$ and $|-r(K+B^*)|$ contracts
$H$ for $r\gg 1$.

Let $h:S'\to S$ be any resolution. The Leray spectral sequence gives an
exact sequence
$$ H^1(S',\o_{S'})\to R^1h_*\o_{S'}\to H^2(S,\o_{S}).
$$
$H^1(S',\o_{S'})=0$ if $S'$ is rational.
$h^2(S,\o_{S})=h^0(S,\o_S(K_S))$ and the latter is zero since
$\o_S(-rK_S)$ has  sections for $r\gg 1$.  Thus $R^1h_*\o_{S'}=0$ and
so $S$ has rational singularities.

Irreducible components of $B^*$ are smooth and rational by the
classification and all other curves in $B^m$ appear as exceptional
curves of $(1,m)$-blow ups, so they are smooth and rational.
\qed

\begin{lem}\label{conn.lem}
 Let $S$ be a real algebraic surface with rational singularities and
$h:S'\to S$ a proper birational morphism from a normal surface $S'$.
Then $S'(\r)\to S(\r)$ has connected fibers.
\end{lem}

Proof. For any $P\in S(\r)$, $H^1(h^{-1}(P), \o)=0$ since $S$ has
rational singularities. Thus $C:=\red(h^{-1}(P))$ is a tree of smooth
rational curves. If $A,B\in C(\r)$ are two points then $C$ has a unique
chain of rational curves  $C_{A,B}$ connecting $A$ and $B$. Complex
conjugation fixes the two ends of the chain, so it fixes every complex
irreducible component and every singular point.  Thus $C_{A,B}(\r)$ is
a connected chain of circles.\qed

\section{The Real Points of Singular Fibers}

In (\ref{S^m-defn.say}) $B^m$ was defined by the formula
$K_{S^m}+B^m\equiv f^* (K_S+B)$. This shows that $S$ is obtained from
$S^m$ by contracting all the curves $A\subset S^m$ such that
$(A\cdot (K_{S^m}+B^m))=0$. By (\ref{S^m-conds}) all such curves are in
$\supp B^m$. 
 Using (\ref{conn.lem}), a version of this also holds for real points.

\begin{prop} \label{top.go.down.prop}
 From $(S^m(\r),B^m(\r))\cong (S^r(\r),B^r(\r))$ we  obtain
$(S(\r),B(\r))$ by     contracting  certain connected subcurves  of
$B^r(\r)$ to   points and by 
   adding   isolated points. \qed
\end{prop}

\begin{rem}\label{top.go.down.rem}
 (1) A real singular point of $S$ may be  isolated in $S(\r)$ and so
invisible in the real part of the resolution.  

(2) In the applications I will be able to compute only a Zariski
neighborhood of $S(\r)$ in $S(\c)$. (\ref{top.go.down.prop}) allows us 
to compute the pair $(S(\r),B(\r))$ up to finite ambiguity.

(3)  The somewhat complicated choice of the partial resolution
$f:S^m\to S$ is important for (\ref{top.go.down.prop}). If we  take the
minimal resolution then the exceptional curves over Du Val points  
appear in $B$ with zero coefficient. Thus
$S^m\to S$ would also involve   exceptional curves which are not
controled by $B^m$.  It is probably possible to analyze this but some
complications definitely do appear.
\end{rem}

We have a list of all possible pairs $(S^r(\r),B^r(\r))$, so using
(\ref{top.go.down.prop}) we can get a list of all possible pairs
$(S(\r),B(\r))$.  From these pairs the real points of the fibers are
obtained by gluing $B$ to itself. It should be kept in mind that if a
point $P\in B$ is glued to its conjugate $\bar P$ then we obtain a real
point. Thus some  points of $X_0(\r)$  may not come from  a point of
$S(\r)$.

The gluing process is easy if
$B(\r)$ is not complicated, say empty or a single circle but it may be 
subtle in general. Some of the worst cases can be avoided if we assume
that $X(\r)$ is orientable though this is not essential.

Next we define a certain class of 2--complexes. The   reason  for this
rather unnatural definition is that this is what I can prove about the
real points of sigular fibers.

\begin{defn} Let $C$ be a circle, $(B_i,\partial B_i)$  2--discs and
$\phi_i:\partial B_i\to  C$ PL--maps. We   use the $\phi_i$ to glue the
discs to $C$. The resulting 2--complex $C\cup_{\phi_i}B_i$ is called a
{\it circle with discs attached}. A disc $B_i$ is called {\it
inessential} if
$\phi_i$ has degree 0.  $\pi_1(C\cup_{\phi_i}B_i)=\z/m$ where $m$ is
the gcd of the degrees of the $\phi_i$. 
\end{defn}

\begin{condition}\label{weird.conds} We consider compact 2--complexes 
$K$ which satisfy the following conditions:
\begin{enumerate}
\item There is a 2--complex $K'$ and a surjective PL--map $h:K'\to K$.
\item $K'$ is the disjoint union of points, intervals, spheres, real
projective planes, tori, Klein bottles and  circles with discs attached.
\item There are finite subsets $A\subset K$ and $A'\subset K'$ such that
$h:K'\setminus A'\to K\setminus A$ is a homeomorphism. Moreover, $A'$
does not contain any interior point of an interval.
\item If $K'$ contains a torus or a Klein bottle then this is its only
2--dimensional connected component.
\end{enumerate} Observe that the 2--dimensional part of $K'$ is
uniquely determined by
$K$. In the applications the lower dimensional parts will play only a
minor role.

It is easy to see that if $K$ satisfies the above conditions and $G$ is
a finite group acting on $K$ without fixed points  then $K/G$ also
satisfies the above conditions. 
\end{condition}

\begin{thm}\label{sing.fib.r.thm}
 Let $X$ be a real algebraic threefold  and $f:X\to C$ a morphism to a
smooth real algebraic curve.  Assume that
 $X$ has isolated singularities which are $\q$-factorial   over $\r$, 
$-K_X$ is
$f$-ample and the smooth part of $X(\r)$ is orientable. Let $0\in
C(\r)$ be a point such that $X_0:=f^{-1}(0)$ 
 is irreducible (over
$\r$). 

Then $X_0(\r)$ is a 2--complex satisfying the conditions
(\ref{weird.conds}). 
\end{thm}

Proof. Assume first that $X_0$ is geometrically reducible. Then $\red
X_0=Z_0+\bar Z_0$ and all real points of $X_0$ are in $Z_0\cap \bar
Z_0$. Thus $X_0(\r)$ is a 1--complex and we are done.

Next consider the case when $X_0$ is geometrically irreducible and
reduced. The list of all   $(S^r(\r),B^r(\r))$ can be read off from
(\ref{S^r-summ.prop}) and  (\ref{top.go.down.prop}) then gives a longer
list for the pairs $(S(\r),B(\r))$. The orientability assumption can be
used through the following two consequences:
\begin{enumerate}
\item  $S^r(\r)\setminus B^r(\r)$ is orientable, and
\item If $\{C_t\}$ is a base point free pencil on $S^r$ such that
$(C_t\cdot (K_{S^r}+B^r))$ is odd then the general $C_t$ intersects at
least one irreducible component of $B^r$ which gets contracted in
$S$.
\end{enumerate} The first of these holds since $S^r(\r)\setminus
B^r(\r)$ injects into
$S(\r)\setminus B(\r)$ by (\ref{top.go.down.prop}) and
$X_0(\r)$ is 2--sided in $X(\r)$, hence its smooth part is orientable.
The second assertion holds since otherwise $\{C_t\}$ would give a base
point free family of curves in $X$ such that 
$(C_t\cdot K_X)$ is odd. This is impossible since $X(\r)$ is orientable
(cf.\ \cite[2.8]{ras}).

A real point of $X_0(\r)$ is either in the image of $S(\r)$ or is in
$\sing(X_0)$. The latter is a 1--complex, and adding a 1--complex does
not change  the conditions (\ref{weird.conds}). 

Let us now go through the list of (\ref{S^r-summ.prop}).

If $S^r$ is a weak Del Pezzo surface such that  the smooth part of
$S^r(\r)$  is orientable then $S^r(\r)$ is obtained from spheres and
tori by identifying some points.  In the smooth case this was proved by
\cite{Comessatti14}, see also \cite{Silhol89}. 
The singular case is discussed in \cite[9.9]{rat3}
or it can be read off   the list of real plane quartics in
\cite{gudkov}. See also \cite{Wall}. 

Next let $S^r$ be a conic bundle. Then $S^r(\r)$ is a union of spheres
with possibly some points identified. $B^r(\r)$ is either empty or a
circle in a sphere. The first case goes as before. If $B^r(\r)$ is
contracted, we get one more sphere.  Otherwise as we go from $S$ to
$X_0$ we have to glue $B^r(\r)$ to itself. Since
$B^r(\r)$ cuts the sphere into 2 discs, we obtain a circle with two
discs attached. If $S^r$ is a conic bundle with one point blown up then
we have spheres and one $\r\p^2$ and $B^r(\r)$ is a line in
$\r\p^2$. If it is contracted then we get a sphere, otherwise we obtain
a circle with one disc attached.

Similar arguments apply every time $B(\r)$ is empty or a circle. Some
of these cases can not happen, for instance $(\r\p^2,\mbox{smooth
conic})$ does not occur since the complement of a smooth conic is not
orientable.

In case  $(\r\p^2,\mbox{two lines})$ a general line in $\p^2$ has
intersection number -1 with $(K_{S^r}+B^r)$. Thus at least one of the
lines gets contracted  in $S$. Similarly, if $S^r=\p^1\times \p^1$ then
all of $B^r$ has to be contracted.

We are left with the case when $S^r$ is a $\p^1$-bundle.  A general
fiber in $S^r$ has intersection number -1 with $(K_{S^r}+B^r)$, and
this shows that the section $H$ gets contracted.  $S^r(\r)$ is a torus
or a Klein bottle. If we contract $H(\r)$ it becomes a sphere with a
pair of points pinched together. If $B^r$ contains at most one real
fiber then we are in one of the previously studied cases.

When $B^r$ contains two real fibers $F_1^r, F_2^r$, further discussions
are needed.  Let $F_i\subset S$ be the birational transform of $F_i^r$. 
Since $H$ is contracted and $B^r$ does not contain any other sections,
$F_1$ and
$F_2$ intersect at a single point $P\in S$.

If $F_1$ and $F_2$ are mapped to the same curve in $X_0$ then we obtain
a circle with 2 discs attached. This is the typical case and it occurs
when $X_0$ is a cone over a nodal rational curve.

We are left with the cases when the normalization $\pi:S\to X_0$ maps
$F_i\to C_i$ and $C_1\neq C_2$. Each $F_i$ appears in $B$ with
coefficient 1, so $F_i\to C_i$ is a degree 2 map and the two branches
of $X$ intersect transversely along $C_i$ by (\ref{what-is-B?}.2).

There are two different degree 2 maps $\r\p^1\to \r\p^1$. The first is
$(z:1)\mapsto (z^2:1)$. Under this map  the circle $\p^1(\r)$ maps to an
interval. If $F_1\to C_1$ is such a map then 
$\pi(S(\r))$ is the circle $C_2(\r)$  with a disc attached and we have
a 1--complex  $C_2(\r)\setminus \pi(F_1(\r))$ added to it.

Finally it may happen that both $F_i(\r)\to C_i(\r)$ are  given as
$(z:1)\mapsto (z:z^2-1)$. I claim that near $\pi(P)$ this leads to a
configuration which can not be realized in $\r^3$.

If $(H^2)$ is odd then  
  locally near $P$
$$ (S(\r), F_1(\r), F_2(\r),P)\sim_s (\r^2, (x=0),(y=0),(0,0)),
$$ where $\sim_s$ denotes a PL homeomorphism of stratified spaces.
  There are  2 more branches of $X_0(\r)$ passing through
$\pi(P)$. One branch $B_1$ passes through $F_1(\r)$ the other
branch $B_2$ passes through $F_2(\r))$ and
these two branches do not intersect each other. Moreover, all
intersections are generically transverse. 
Intersect everything with a small sphere
$S^2_{\epsilon}$
around the origin. The upper hemi sphere $(z\geq 0)$
is a topological disc and $B_1$ and $B_2$ intersect it in 2 curves.
One of these curves connects the points $(0,1),(0,-1)$
and the other the points $(1,0),(-1,0)$. Hence 
$B_1\cap B_2\cap S^2_{\epsilon}\neq \emptyset$, a contradiction.

If $(H^2)$ is even then   locally near $P$
\begin{eqnarray*}
\lefteqn{(S(\r), F_1(\r), F_2(\r),P)\sim_s}\\ && \quad\sim_s
((xy-z^2=0), (x=z=0), (y=z=0),(0,0,0))\subset \r^3.
\end{eqnarray*}
 There are  2 more branches of
$X_0(\r)$ passing through
$\pi(P)$. One passes through $F_1(\r)$ the other through $F_2(\r)$ and
these two branches do not intersect each other. Moreover, all
intersections are generically transverse. This is again impossible in
$\r^3$. 
\medskip

We are left with the case when $X_0$ is a multiple fiber. Let $m>0$ be
the smallest integer such that $mX_0$ is Cartier. Then $\o_X(mX_0)$ is
a locally free sheaf isomorphic to $\o_X$ in a neighborhood of $X_0$. 
The section $1\in H^0(X,\o_X)\cong H^0(X,\o_X(mX_0))$ determines the
corresponding $m$-sheeted cyclic cover
 $Y\to X$   (cf. \cite[Sec.\ 2.4]{KoMo98}). Since $X_0$ is Cartier at
real points,  $Y\to X$ is unramified at real points. Thus 
$Y(\r)\to X(\r)$ is a homeomorphism if $m$ is odd. If $m$ is even then
over each connected component we get either a 2--sheeted cover or the
empty set. In the latter case we switch to $-1\in H^0(X,\o_X)$.

We have already proved that $Y_0(\r)$ satisfies the properties
(\ref{weird.conds}), thus the same holds for $X_0(\r)$, as remarked at
the end of (\ref{weird.conds}).
\qed

\section{Proof of the Main Theorem}

In this section I use (\ref{sing.fib.r.thm}) and a purely topological
argument to complete the proof of (\ref{main.dp.thm}). First, following
(\ref{2nd.red.say}),  we identify the neighborhoods of the singular
fibers. Then we show how these pieces are assembled to form
$\overline{X(\r)}$.

\begin{lem}\label{main.thm.pf.1}
 Let $M$ be an orientable 3--manifold and $K\subset M$ a 2-complex
satisfying the conditions (\ref{weird.conds}). Let
$K\subset U\subset M$ be a regular neighborhood of $K$ and assume that
$\partial U$ is a union of spheres and tori. Then every connected
component of $U$ is one of the following.
\begin{enumerate}
\item connected sum of lens spaces and $S^1\times S^2$ minus balls,
\item connected sum of lens spaces,   $S^1\times S^2$  and of a solid
torus minus balls,
\item  interval bundle over a torus or a Klein bottle.
\end{enumerate}
\end{lem}

Proof. As a first step we replace $h:K'\to K$ with another map
$\tilde h:\tilde K'\to \tilde  K$ such that 
\begin{enumerate}\setcounter{enumi}{3}
\item $\tilde h:\tilde K'\to \tilde  K$  satisfies the conditions 
(\ref{weird.conds}),
\item $U$ is a regular neighborhood of $\tilde  K$,
\item $\tilde h$ is injective on the set of 2--dimensional points of
$\tilde  K'$.
\end{enumerate}

To achieve this, note that there are only finitely many points $p\in K$
over which $h$ is not one--to--one. We will get rid of these one at a
time. Let $p\in V$ be a regular neighborhood. $K\cap \partial V$ is a
union of connected 1--complexes $A_j$ and $K\cap V$ is the cone over
$K\cap \partial V$. If all the $A_j$ have dimension 0 then we do not
need to do anything. Otherwise, there is  a 1--dimensional component
(say $A_1$) and a PL-homeomorphism 
$\psi:(V,\partial V)\to (B^3,S^2)$   such that $\psi(A_1)$ is in the
northern hemisphere and all the other 1--dimensional $A_j$ map to the
southern hemisphere. Now   cut the 3--ball $B^3$ along the equator and
move the two halves apart a little. The center of the ball  sweeps out a
small interval. We add this interval to obtain $K_1$. 
$K'_1$ is $K'$ union an interval.

Repeating this procedure if necessary, at the end we obtain 
$\tilde h:\tilde K'\to \tilde  K$. 

To simplify notation let us assume that $h:K'\to K$ already satisfies
the above conditions 4--6. Next we use induction on the number of
2--dimensional components of $K'$.

If there is an $S^2\sim L\subset K'$ then $h(L)$ is an embedded $S^2$.
Cut $U$ along $h(L)$ and glue 3--balls to the resulting two spheres to
get
$U_1$. 
$K'$ is replaced by $K'_1:=K'\setminus L$. $U$ is obtained from $U_1$ by
attaching a 1--handle which is either taking connected sum of two
components or taking connected sum with $S^1\times S^2$.  If there is
an 
$\r\p^2\sim L\subset K'$ then
$h(L)$ is an embedded
$\r\p^2$ and the boundary of its regular neighborhood is $S^2$ (since
$M$ is orientable). Cutting along this $S^2$ corresponds to connected
sum with
$\r\p^3$ (which is the lens space $L_{2,1}$). 

Assume now that there is an $L \subset K'$ which is a circle with discs
attached. Let $V$ be a regular neighborhood of $h(C)$ which we may
assume to be an embedded circle.
$V$ is a solid torus.  Consider the case when one of the discs  (say
$B_1$) is inessential. Then $\partial V\cap h(N_1)$ bounds a  disc 
$D\subset
\partial V$. We may asume that no other $B_i$ intersects this disc. We
change $h(B_1)$ by replacing $V\cap B_1$ with $D$ and adding an interval
connecting $D$ with $h(C)$ to $K$. This way we have not changed $U$, we
have removed one inessential disc and we created a new embedded sphere.
By repeating this procedure we may assume that there are no inessential
discs. Then each $h(B_i)$ intersects $\partial V$ in a simple closed
curve which is not null homotopic. Since these curves are all disjoint,
they are in the same homotopy class $\gamma$. Thus  $V\cup_ih(B_i)$ is
a solid torus with discs attached along parallel curves in $\partial
V$.  The boundary of a regular neighborhood $W$ of $V\cup_ih(B_i)$ is a
union of spheres. We can again cut $U$ along these spheres.  $W$ can
also be obtained as attaching first a solid torus $V'$ to $V$  such
that the meridian of $V'$ maps to $\gamma$ and then removing some balls
from $V'$. Gluing two solid tori gives a lens space (\ref{lens.defn}).

We are left to deal with the 1--dimensional part of $K$. Every connected
component can be collapsed to a bouquet of circles. If there are $\geq
2$ circles then the boundary of its regular neighborhood has genus
$\geq 2$, a contradiction.  Thus we get either  a ball (if ther are no
cicrcles) or a solid torus (if there is one circle).

If there is an $L\subset K'$ which is a torus or a Klein bottle then
the boundary of a regular neighborhood of $h(L)$ is either 2 or 1 tori.
There are no  other 2--dimensional components in $K'$ by
(\ref{weird.conds}.4).  If $K$ has a 1--dimensional subcomplex which
does not collapse into $h(L)$ then it leads to a genus $\geq 2$
component in
$\partial U$, a contradiction. \qed

\begin{rem} Using (\ref{what-is-B?}) we obtain that the only lens spaces
that can appear in  (\ref{main.thm.pf.1}) are $\r\p^3=S^3/\z_2$ and
$S^3/\z_3$. I have no example for the latter.
\end{rem}

\begin{say}[Proof of (\ref{main.dp.thm})]\label{main.thm.pf.2}{\ }

We have established that $M$ is glued together from the following
pieces:
\begin{enumerate}
\item $S^1\times S^2$ minus open balls,
\item  lens space minus open balls,
\item solid torus minus open balls,
\item interval bundle over a torus or a Klein bottle.
\end{enumerate}

An orientable  interval bundle over a torus is a torus times an
interval and   the boundary of an orientable  interval bundle over a
Klein bottle is  a torus (cf.\ \cite[1.6]{rat2}).

If the gluing involves  a torus times an interval then this piece can be
thrown away, except when the two boundary components are glued to each
other. This gives a torus bundle over  a circle.

Two copies of an  interval bundle over a  Klein bottle glued together
map to an interval such that the fibers are tori over interior porints
and Klein bottles over the two boundary points. This 3--manifold is
doubly covered by a torus bundle over  a circle.

An interval bundle over a  Klein bottle glued to  a solid torus is
homeomorphic to $\r\p^3\# \r\p^3$ (cf.\ \cite[12.7]{rat2}), so in this
case we can change the decomposition to one that does not involve   any
interval bundles over a  Klein bottle. 

We are left with the case when $M$ is glued together from  lens spaces
minus open balls.  We do one gluing at a time. If a new lens space is
glued in, that is connected sum. If two boundary components of a
connected component are glued together then that is the same as taking
connected sum with
$S^1\times S^2$. \qed
\end{say}

\section{Examples}

\begin{exmp}\label{deg2.deg} Consider $\p^4(1,1,1,2,2)$ with
coordinates $(x:y:z:u:v)$ and the affine line $\a^1$ with coordinate
$t$.  Let $X$ be the complete intersection
$$
\begin{array}{l} X\subset \p^4(1,1,1,2,2)\times \a^1
\qtq{given by equations} \\ u^2+v^2=f_4(x,y,z)\qtq{and} tv=q_2(x,y,z).
\end{array}
$$ Let $\pi:X\to \a^1$ be the second projection. For $t\neq 0$ we can
eliminate $v$ to obtain  a degree 2 Del Pezzo surface
$$ X_t\cong (u^2=f_4(x,y,z)-t^{-2}q_2(x,y,z)^2)\subset \p^3(1,1,1,2).
$$ For $t=0$ the equations become $u^2+v^2=f_4(x,y,z)$ and
$q_2(x,y,z)=0$. This has  two points of index 2 at
$(0:0:0:1:\pm\sqrt{-1})$, both are analytically isomorphic to
$\c^2/\z_4(1,1)$. The projection 
$(x:y:z:u:v)\mapsto (x:y:z)$ is defined outside these two points and
the minimal resolution of $X_0$ becomes a conic bundle over the conic
$(q_2(x,y,z)=0)$. The singular fibers correspond to the solutions of
the equations $q_2=f_4=0$. We get various cases depending on how the
curves $(q_2=0)$ and $(f_4=0)$ intersect.

Another way to obtain this model is as follows. Consider the  family
$$ Y\subset \p^3(1,1,1,2)\times \a^1
\qtq{given by}   u^2=t^2f_4(x,y,z)-q_2(x,y,z)^2.
$$ Outside the origin $Y$ is isomorphic to $X$ via the transformation
$$ (x,y,z,u,t)\mapsto (x,y,z,ut,t).
$$ The central fiber  $Y_0$ consists of a conjugate pair of planes
intersecting along the conic $(q_2=0)$.
$Y$ is singular along $(u=t=q_2=0)$. If we blow up the singular curve,
the two planes in the central fiber become disjoint. Contracting them
 gives the 3--fold $X$. 
\end{exmp}

\begin{exmp}\label{deg1.deg} Consider $\p^4(1,1,2,3,3)$ with
coordinates $(x:y:z:u:v)$ and the affine line $\a^1$ with coordinate
$t$.  Let $X$ be the complete intersection
$$
\begin{array}{l} X\subset \p^4(1,1,2,3,3)\times \a^1
\qtq{given by equations} \\ u^2+v^2=f_6(x,y,z)\qtq{and} tv=c_3(x,y,z).
\end{array}
$$ Let $\pi:X\to \a^1$ be the second projection. For $t\neq 0$ we can
eliminate $v$ to see that the fiber  is the degree 1 Del Pezzo surface
$$ X_t\cong (u^2=f_6(x,y,z)-t^{-2}c_3(x,y,z)^2)\subset \p^3(1,1,2,3).
$$ For $t=0$ the equations become $u^2+v^2=f_6(x,y,z)$ and
$c_3(x,y,z)=0$. This has  two points of index 3 at
$(0:0:0:1:\pm\sqrt{-1})$, both are analytically isomorphic to
$\c^2/\z_9(1,2)$. The projection 
$(x:y:z:u:v)\mapsto (x:y:z)$ is defined outside these two points and
the minimal resolution of $X_0$ becomes a conic bundle over the curve
$(c_3(x,y,z)=0)$ blown up in one point. Although $c_3$ has degree 3,
this curve is birational to $\p^1$. Indeed, since $z$ has degree 2, it
appears  in $c_3$ only linearly, so $z$ can be rationally expressed in
terms of
$x,y$. 

One can also  obtain this model from the degeneration
$$ Y\subset \p^3(1,1,2,3)\times \a^1
\qtq{given by}   u^2=t^2f_6(x,y,z)-c_3(x,y,z)^2,
$$ but the birational transformation between them is   more complicated.
\end{exmp}

\begin{exmp} Start with the trivial family of quadrics
$$
\p^1\times \p^1\times \a^1\qtq{with coordinates}  (x_1:x_2),(y_1:y_2),
t.
$$ Consider the $\z_2$-action
$$
\tau_1: (x_1,x_2,y_1,y_2,t)\mapsto  (x_2,-x_1,y_2,-y_1,-t).
$$ This has 4  fixed points  at $(1,\pm\sqrt{-1},1,\pm\sqrt{-1},0)$.
Set $X_1:=(\p^1\times \p^1\times \a^1)/\z_2(\tau_1)$.
$X_1$ has 4 singularities of analytic type $\c^3/\z_2(1,1,1)$. The
central fiber is double and the reduced central fiber has 4
$A_1$-type points. It is easy to see that $X_1(\r)$ is not orientable
and its central fiber is $S^1\times S^1$.

Another $\z_2$-action is given by
$$
\tau_2: (x_1,x_2,y_1,y_2,t)\mapsto  (x_2,-x_1,y_1,-y_2,-t).
$$ This also has 4  fixed points   at $(1,\pm\sqrt{-1},1,0,0)$ and
$(1,\pm\sqrt{-1},0,1,0)$. Set $X_2:=(\p^1\times \p^1\times
\a^1)/\z_2(\tau_2)$.
$X_2$ has 4 singularities of analytic type $\c^3/\z_2(1,1,1)$, the
central fiber is double and the reduced central fiber has 4
$A_1$-type points. One can see that $X_2(\r)$ is   orientable and  its
central fiber is a Klein bottle.
\end{exmp}

\begin{exmp}  Let $(x:y:z)$ and $t$ be coordinates on $\p^2\times \a^1$
and consider the surface
$S:=(x^2+y^2-tz^2=0)$. Let $\z_n$ act on $(x,y,z,t)$ as a rotation  of
order $n$ on $x,y$ and  identity on $z,t$. This induces an action on
$\o_S+\o_S$. Let $F\subset \o_S+\o_S$ be the locally free rank 2
subsheaf which on the $z\neq 0$ affine chart is generated by the
sections $(x,y)$ and $(-y,x)$.  Outside $(x^2+y^2=0)$ the subsheaf $F$
is the same as $\o_S+\o_S$. Let $X:=\p_SF$ be the corresponding
$\p^1$-bundle over $S$. 
 $\z_n$ acts on $X$. This action has two isolated conjugate fixed
points over $(0,0,1,0)\in S$ and fixes the conjugate surfaces over the
curves  $(1,\pm\sqrt{-1},0,t)\subset S$.  Set $X_n:=X/\z_n$.  The
projection $\pi: X_n\to \a^1$ exhibits $X_n$ as a degeneration of
quadrics and $-K_{X_n}$ is $\pi$-ample.
$X_n$ has only terminal singularities, it has 2 points of index $n$ 
and it is smooth at all real points.  The central fiber is
geometrically reducible.

The only slight problem with this example is that $\rho(X_n/\a^1)=2$
since $\pi$ can be factored as $X_n\to S/\z_n\to \a^1$.  Probably one
can  globalize this example  to get  relative Picard number 1.
\end{exmp}

\begin{exmp} Let $S^*\subset \p^3$ be a smooth quadric and $B^*=0$. Let
$P,\bar P\in S^*$ be a pair of conjugate points not on a line. Let
$S_1\to S^*$ be the blow up of $P+\bar P$ with exceptional curve
$E+\bar E$. Set $B_1:=(1/2)(E+\bar E)$. The pencil of planes through
$P+\bar P$ gives $S_1$ a conic bundle structure with $E+\bar E$ as
conjugate sections. Blow up 3 more pairs of conjugate points on $E+\bar
E$ to get $S^m$.  The birational transforms of $E,\bar E$ have self
intersection $-4$, so they can be contracted
 $S^m\to S$. This $S$ is among the types described in (\ref{deg2.deg}).
\end{exmp}

\begin{exmp} Set  $Q:=\p^1_{(x_0:x_1)}\times \p^1_{(y_0:y_1)}$ and let 
$\sigma:Q\to Q$ be given by $(x_0:x_1, y_0:y_1)\mapsto
(y_1:y_0,x_0:x_1)$. Then $\sigma$ has order 4 on $Q$ and   on
$H_1(Q(\r),\z)$. Take $Q\times  \p^1_{(z_0:z_1)}$. Let $\tau$ be the
action which is
$\sigma$ on $Q$ and 
$(z_0:z_1)\mapsto (z_0+z_1: -z_0+z_1)$ on the second factor. Set
$X:=(Q\times  \p^1_{(z_0:z_1)})/(\tau)$. Then
$X\to \p^1_{(z_0:z_1)}/(\tau)$ is a quadric bundle (except over a
conjugate pair of complex points). Its real part gives  an $S^1\times
S^1$-bundle over $S^1$ with order 4 monodromy.
\end{exmp}

\begin{exmp} Let $Q$ be  the degree 6 Del Pezzo surface given in 
$\p^1_{(x_0:x_1)}\times
\p^1_{(y_0:y_1)}\times \p^1_{(z_0:z_1)}$ by the equation
$$ x_0y_0z_1+x_0y_1z_0+x_1y_0z_0= x_0y_1z_1+x_1y_0z_1+x_1y_1z_0.
$$ Let 
$\sigma:Q\to Q$ be given by 
$$ (x_0:x_1, y_0:y_1,z_0:z_1)\mapsto (z_1:z_0,x_1:x_0,y_1:y_0).
$$
 $\sigma$ has order 6 on $Q$ and   on $H_1(Q(\r),\z)$. Take $Q\times 
\p^1_{(t_0:t_1)}$. Let $\tau$ be the action which is
$\sigma$ on $Q$ and 
$$ (t_0:t_1)\mapsto (\cos(\pi/6)t_0+\sin(\pi/6)t_1:
-\sin(\pi/6)t_0+\cos(\pi/6)t_1)
$$
 on the second factor. Set
$X:=(Q\times  \p^1_{(t_0:t_1)})/(\tau)$. Then
$X\to \p^1_{(t_0:t_1)}/(\tau)$ is a  degree 6 Del Pezzo surface bundle
(except over a conjugate pair of complex points). Its real part gives 
an $S^1\times S^1$-bundle over $S^1$ with order 6 monodromy.
\end{exmp}

\begin{exmp} In $\p^3_{(x_0:x_1:x_2:x_3)}\times
\p^1_{(y_0:y_1)}$ consider the 3--fold
$$ X:=(y_0(x_0^2+x_1^2)=y_1(x_2^2+x_3^2)).
$$ For $p$ odd  and $(p,q)=1$ let 
$\sigma:X\to X$ be the action which is rotation by
$2\pi/n$ on $(x_0,x_1)$ and  rotation by
$2q\pi/n$ on $(x_2,x_3)$. Then $X(\r)\sim S^3$ and $X(\r)/(\sigma)\sim
L_{p,q}$. 

The complex 3--fold $X/(\sigma)$  is smooth at its real points but it 
has  nonterminal singularities at complex points. At least for $q=1$
these are easy to resolve and one obtains a quadric bundle $X'$ such
that $X'(\r) \sim L_{p,1}$.  The other cases seem more complicated.
\end{exmp}

\noindent University of Utah, Salt Lake City UT 84112 

\begin{verbatim}kollar@math.utah.edu\end{verbatim}

\end{document}